\documentclass[10pt,leqno]{amsart}
\usepackage{amscd}
\usepackage{amssymb}
\usepackage{amsfonts}
\usepackage{latexsym}
\usepackage{verbatim}

\theoremstyle{plain}
\newtheorem{theorem}{Theorem}[section]
\newtheorem{definition}[theorem]{Definition}
\newtheorem{lemma}[theorem]{Lemma}
\newtheorem{prop}[theorem]{Proposition}

\newtheorem{rem}[theorem]{Remark}
\newtheorem{ex}[theorem]{Example}
\begin{document}
\title{On Symplectc half-flat manifolds}
\author{Adriano Tomassini and Luigi Vezzoni}
\date{\today}
\address{Dipartimento di Matematica\\ Universit\`a di Parma\\ Viale
  G. P. Usberti 53/A\\
43100 Parma\\ Italy}
\email{adriano.tomassini@unipr.it}
\address{Dipartimento di Matematica \\ Universit\`a di Torino\\
Via Carlo Alberto  10\\
10123 Torino, Italy}
\email{luigi.vezzoni@unito.it}
\thanks{This work was supported by the Projects M.I.U.R. ``Geometric Properties of Real and Complex Manifolds'',
``Riemann Metrics and  Differenziable Manifolds'' and by G.N.S.A.G.A.
of I.N.d.A.M.}
\begin{abstract}
We construct examples of symplectic half-flat manifolds on compact quotients of solvable Lie groups. We prove that the Calabi-Yau
structures are not rigid in the class of symplectic half-flat structures. Moreover, we provide an example of a  compact $6$-dimensional
symplectic half-flat manifold whose real part of the complex volume form is $d$-exact. Finally we discuss the $4$-dimensional case.
\end{abstract}
\maketitle
\newcommand\C{{\mathbb C}}
\newcommand\R{{\mathbb R}}
\newcommand\Z{{\mathbb Z}}
\newcommand\T{{\mathbb T}}
\newcommand{\de}[2]{\frac{\partial #1}{\partial #2}}
\newcommand\Span{{\rm Span}}
\section{Introduction}
\noindent Half-flat structures arise as a special class of
SU(3)-structures introduced by Hitchin in \cite{Hi1}. An
SU$(3)$-structure $(\omega,J,\psi)$ on a $6$-dimensional manifold
$M$ is said to be \emph{half-flat} if the defining forms
$\omega\in\Lambda^2(M)$, $\Re\mathfrak{e}\,\psi\in\Lambda^3(M)$
satisfy
\begin{equation}\label{HL}
d\omega\wedge\omega=0\,,\quad d\Re\mathfrak{e}\,\psi=0\,.
\end{equation}
Condition \eqref{HL} is equivalent to require that the intrinsic torsion of $(\omega,J,\psi)$ is symmetric.
In \cite{Hi1} Hitchin proves that, starting with a
half-flat manifold $(M,\omega,J,\psi)$, if certain evolution equations
have a solution coinciding with the initial datum $(\omega,\psi)$ at
time $t=0$, then there exists a metric with holonomy contained in
$G_2$ on $M\times I$ for some interval $I$ (see also \cite{CS},
\cite{CF}).\\

In the present paper we study {\em symplectic half-flat manifolds}, namely $6$-dimensional
manifolds $M$ endowed with an SU(3)-structure $(\omega,J,\psi)$ satisfying the following
\footnote{In \cite{dBT2,dBT3} symplectic half-flat manifolds were called
\emph{special generalized Calabi-Yau manifolds}.
Here we change the terminology to avoid confusion with the case considered by Hitchin in \cite{Hi3}.}
\begin{equation}\label{specialsymplectic}
d\omega =0\,,\quad d\Re\mathfrak{e}\,\psi=0\,,\quad \psi\wedge\overline{\psi}=-\frac{4}{3}i\omega^3\,.
\end{equation}
In such a case, since $\Re\mathfrak{e}\,\psi$
is a calibration on $M$, we can define {\em special Lagrangian
submanifolds} as compact three-dimensional Lagrangian submanifolds
of $M$, calibrated by $\Re\mathfrak{e}\,\psi$
(see \cite{HL}). \\
It turns out (see lemma 2.3 and \cite{dBT2} for its proof) that the
complex volume form $\psi$ of a symplectic half-flat manifold
$(M,\omega,J,\psi)$ is parallel with respect to the Chern connection
$\nabla$ of the almost K\"ahler manifold $(M,\omega,J)$. Therefore,
in such a case, the holonomy of $\nabla$ is contained in SU$(3)$.
Note that if the almost complex structure $J$ is integrable, then
$(M,\omega,J)$ is a K\"ahler manifold and the Chern connection
coincides with the Levi-Civita one and, consequently, in this case symplectic
half-flat manifolds are Calabi-Yau manifolds of complex dimension
$3$. Therefore, symplectic half-flat manifolds can be viewed as an
extension of Calabi-Yau manifolds to the non-integrable case. In
particular, symplectic half-flat manifolds are {\em generalized
Calabi-Yau manifolds} in the sense of Hitchin. Indeed, Hitchin in
\cite{Hi3} gives the notion of {\em generalized Calabi-Yau
structure} on a manifold $M$ of dimension $2m$ as a closed complex
form $\varphi$ of mixed degree which is a complex pure spinor for
the orthogonal vector bundle $TM\oplus T^*M$ endowed with the
natural pairing $<,>$ and such that $<\varphi
,\overline{\varphi}>\neq 0$. It turns out that such a structure (see
\cite{Hi3}) induces a {\em generalized complex structure} on $M$
(see \cite{G}). Two basic examples of generalized Calabi-Yau
manifolds are furnished by holomorphic Calabi-Yau manifolds and
symplectic manifolds. In the first case if $\Psi$ denotes the
holomorphic volume form on the K\"ahler manifold $M$, then $\varphi
=\Psi$ defines a generalized Calabi-Yau structure on $M$. In the
second case, if $\omega$ is the symplectic form on $M$, then
$\varphi =\exp(i\omega)$ is a generalized Calabi-Yau structure on
$M$.

In this paper we give some explicit examples of compact
symplectic half-flat manifolds arising as quotient of solvable Lie groups and
we construct a family of symplectic half-flat structures $(\omega_t,J_t,\psi_t)$ on the $6$-dimensional
torus $\T^6$, coinciding with the standard Calabi-Yau structure for $t=0$, but which is not Calabi-Yau for $t\neq 0$ (see example \ref{torus}
and theorem \ref{remtorus}). For other examples we refer to \cite{dBT3} and \cite{CT}.
Then we give a symplectic half-flat structure on a compact $6$-manifold whose complex volume form
has real part $d$-exact (see theorem \ref{principale}). This is in contrast with the integrable case, namely in the context of Calabi-Yau
geometry, where the real part of the complex volume form cannot be exact.
Finally, we consider the $4$-dimensional case. As in the $6$-dimensional case,
we characterize such structures in terms of stable forms (see proposition \ref{SU2}).\\
In section 2, we start by recalling some facts on symplectic half-flat manifolds. In section 3
we describe the compact examples previous mentioned. In the last section we consider the 4-dimensional case.
\section{Symplectic half-flat geometry} 
Let $M$ be a $6$-dimensional smooth manifold and let $L(M)$ be the principal $\mbox{GL}(6,\R)$-bundle of
linear frames on $M$. An $\mbox{SU}(3)$-structure on $M$ is a reduction of $L(M)$ to a principal bundle whose structure group
is isomorphic to $\mbox{SU}(3)$. It is well known that $\mbox{SU}(3)$-structures on $M$ are in one-to-one correspondence with the triple
$(\omega,J,\psi)$, where
\begin{itemize}
\item $\omega$ is a non-degenerate $2$-form on $M$;
\vskip0.1cm
\item $J$ is an $\omega$-calibrated almost complex structure, i.e. $J$ is an almost complex structure on $M$ such that the tensor
$
g_J(\cdot,\cdot)=\omega(\cdot,J\cdot)
$
is an almost Hermitian metric.
\vskip0.1cm
\item $\psi\in\Lambda^{3,0}_J(M)$ is a complex volume form on $M$ such that
$$
\psi\wedge\overline{\psi}=-\frac{4}{3}i\omega^3\,.
$$
\end{itemize}
\begin{definition}
An $\mbox{\emph{SU}}(3)$-structure $(\omega,J,\psi)$, is said to be \emph{symplectic half-flat} if
$$
d\omega=0\,,\quad d\Re\mathfrak{e}\,\psi=0\,.
$$
\end{definition}
Symplectic half-flat manifolds lie in the intersection of symplectic manifolds and half-flat manifolds. The latter ones have been introduced
by Hitchin in \cite{Hi1} (see also \cite{CS}).


Now we present symplectic half-flat structures in terms of differential forms and the Chern connection.
In order to do this we start with recalling the following
\begin{definition}
Let $(V,\omega)$ be a symplectic vector space of dimension $2n$; the
{\em symplectic Hodge operator}
$$
\bigstar\colon \Lambda^r(V^*)\to\Lambda^{2n-r}(V^*)
$$
is defined by
$$
\alpha\wedge\bigstar\beta=\omega(\alpha,\beta)\frac{\omega^n}{n!}\;.
$$
\end{definition}
\noindent Now we describe the standard model:\\
On $\R^6$ let us consider the natural symplectic structure
$$
\omega_3=dx_1\wedge dx_4+dx_2\wedge dx_5+dx_3\wedge dx_6\,,
$$
where $\{x_1,\dots,x_6\}$ are the standard coordinates on $\R^6$, and let
$$
\Lambda_0^3(\R^{6*})=
\{\Omega\in\Lambda^3(\R^{6*})\;\;|\;\;\Omega\wedge\omega_3=0\}\,.
$$
The group
$$
G=\hbox{\rm Sp}(3,\R)\times\R^*_+
$$
acts on $\Lambda^3_0(\R^{6*})$ by
$$
(A,t)(\Omega)=tA^*(\Omega)\,,
$$
for any $A\in \hbox{\rm Sp}(3,\R)$, $t\in\R_+^*$, where $\hbox{\rm Sp}(3,\R)$ denotes the symplectic group
on $(\R^6,\omega_3)$.\\
Moreover, the 3-form
$$
\Omega_0=\Re \mathfrak{e}\,(dz_1\wedge dz_2\wedge dz_3)
$$
belongs to $\Lambda^3_0(\R^{6*})$,
where $dz_h=dx_h+idx_{h+3}$, $h=1,2,3$.\\
Let us consider the map
$$
\Lambda_0^3(\R^{6*})\to \hbox{\rm End} (\Lambda^1(\R^{6*}))
$$
defined by
$$
P_{\Omega}(\alpha)=-\frac{1}{2}\bigstar(\Omega\wedge\bigstar
(\Omega\wedge\alpha))\,,
$$
for any $\alpha\in\Lambda^1(\R^{6*})$. We have
\begin{itemize}
\item $\omega_3(P_{\Omega}\alpha,\beta)=-\omega
_3(\alpha,P_{\Omega}\beta)$,
\vspace{0.2cm}
\item $P_{\Omega}^2=c\,I$, $c\in\R$
\end{itemize}
(see e.g. \cite{dBT2}).\\
Given $\Omega\in\Lambda^3_0(\R^{6*})$, let
$$
\begin{aligned}
F_{\Omega}\colon\Lambda^1(\R^{6*})&\to\Lambda^4(\R^{6*})\\
F_{\Omega}(\alpha)&=\Omega\wedge\alpha.
\end{aligned}
$$
Then we get that the following facts are equivalent:
\begin{itemize}
\item $\Omega$ belongs to the $G$-orbit of $\Omega_0$,
\vspace{0.2cm}
\item $F_{\Omega}$ is injective and $\omega_3$ is negative defined
  on $\hbox{\rm Im}(F_{\Omega})$.
\end{itemize}
A 3-form $\Omega\in \Lambda^3(\R^{6*})$ is said to be {\em positive} if $F_{\Omega}$ is injective and $\omega_3$ is negative
defined on $\hbox{\rm Im}(F_{\Omega})$
(in this case $c =-\sqrt[3]{\det P_{\Omega}}$) and \emph{normalized} if
$\det P_{\Omega}=1$.

Let $(M,\omega)$ be a symplectic manifold. Let $J$ be an $\omega$-calibrated almost complex structure on $M$ and let
$$
g_J(\cdot,\cdot)=\omega(\cdot,J\cdot)
$$
be the almost Hermitian metric associated with $(\omega,J)$. Denote by $\nabla^{LC}$ the
Levi-Civita connection of $g_J$. Then the {\em Chern connection} on
$(M,\omega,J)$ is defined by
$$
\nabla=\nabla^{LC}-\frac{1}{2}J\nabla^{LC}J\,.
$$
It is known that
$$
\nabla g=0\,,\quad\nabla J=0\,,\quad T^{\nabla}=\frac{1}{4}N_J
$$
where
$$
N_J(X,Y)=[JX,JY]-J[JX,Y]-J[X,JY]-[X,Y]
$$
is the {\em Nijenhuis tensor of $J$} (see e.g. \cite{Gau}). Let $(M,J)$ be an almost
complex manifold; then the exterior derivative
$$
d\colon\Lambda_J^{p,q}(M)\to\Lambda_J^{p+2,q-1}(M)\oplus\Lambda_J^{p+1,q}
(M)\oplus\Lambda_J^{p,q+1}(M)\oplus\Lambda_J^{p-1,q+2}(M)
$$
splits as
$$
d=A_J+\partial_J+\overline{\partial}_J+\overline{A}_J\,.
$$
Note that by the Newlander-Nirenberg theorem $J$ is complex if
and only if $A_J=0$.

Let $(M,\omega)$ be a $6$-dimensional (compact) manifold  equipped with an
$\omega$-calibrated almost complex structure
$J$ and a complex $(3,0)$-form $\psi$ satisfying
$$
\psi\wedge\overline{\psi}=-ie^\sigma\,\omega^3\,,
$$
where $\sigma$ is a $C^\infty$ function on $M$. Set
$$
\Omega :=\Re \mathfrak{e}\,\psi\,.
$$
The following lemma (see \cite{dBT2})
gives a characterization of symplectic half-flat structures in terms of differential forms and the Chern connection.
\begin{lemma}
\label{lammadBT1}
With the notation above, the following facts are equivalent
\begin{enumerate}
\item[a)] $\begin{cases}
d\Omega=0\\
\sigma=const.
\end{cases}$\vskip.2truecm\noindent
\item[b)]
$\begin{cases}
d\Omega =0\\
\Omega\wedge\omega = 0 \\
\frac{2}{3}\sqrt{3}\,e^{-\sigma/2}\Omega \hbox{ is positive and normalized at any point}
\end{cases}$\vskip.2truecm\noindent
\item[c)]
$\begin{cases}
\nabla\psi=0\\
A_J(\overline{\psi})+\overline{A}_J(\psi)=0\,.
\end{cases}$
\end{enumerate}
\end{lemma}
\begin{rem}{\rm By lemma \ref{lammadBT1} one can define a
symplectic half-flat manifold as a (compact) $6$-dimensional
symplectic manifold $(M,\omega)$ endowed with an $\omega$-calibrated
almost complex structure $J$ and a (3,0)-form $\psi$ such that
$$
\begin{cases}
\nabla\psi=0\\
A_J(\overline{\psi})+\overline{A}_J(\psi)=0\,.
\end{cases}
$$
Therefore, the holonomy of the Chern connection of a symplectic half-flat manifold is
contained in SU$(3)$.

Alternatively, a symplectic half-flat structure can be given by a pair $(\omega,\Omega)$, where
$\omega$ is a symplectic form and $\Omega$ is a closed 3-form satisfying $\Omega\wedge\omega=0$ and which
is positive and normalized at any point. Indeed, in this case the $\omega$-calibrated
almost complex structure is given by the isomorphism dual to $P_\Omega$ and the complex volume form is
$$
\psi:=\Omega+iP_\Omega(\Omega)\,.
$$
}
\end{rem}
\begin{rem}
{\rm Observe that if
$$
d\Re \mathfrak{e}\,\psi=d\Im \mathfrak{m}\,\psi=0\,,
$$
then the almost complex structure is integrable (hence
$(M,\omega,J,\psi)$ is a Calabi-Yau manifold).\\
Indeed, if $\alpha\in\Lambda^{1,0}_J(M)$, then
$$
0=d(\alpha\wedge\psi)=d\alpha\wedge\psi
$$
and consequently
$$
d(\Lambda^{1,0}_J(M))\subset\Lambda^{2,0}_J(M)\oplus\Lambda^{1,1}_J(M).
$$
}
\end{rem}
\vskip0.3cm Let $(M,\omega,J,\psi)$ be a symplectic
half-flat manifold. As in the Calabi-Yau case, the 3-form
$\Omega=\Re\mathfrak{e}\,\psi$ is a calibration on $M$ (see
\cite{HL}). Therefore we can give the following
\begin{definition}
A {\em special Lagrangian submanifold} of $(M,\omega,J,\psi)$ is
a compact submanifold $p\colon L\hookrightarrow M$ calibrated by
$\Re\mathfrak{e}\,\psi$.
\end{definition}
\noindent We have the following
\begin{lemma}
Let $p\colon L\hookrightarrow M$ be a submanifold. The following
facts are equivalent
\begin{enumerate}
\item[1.] $p^*(\omega)=0$, $p^*(\Im\mathfrak{m}\,\psi)=0$;
\vskip0.2cm
\item[2.] there exists an orientation on $L$ making it calibrated by
$\Re\mathfrak{e}\,\psi$.
\end{enumerate}
\end{lemma}
\section{Examples of compact symplectic half-flat solvmanifolds}
\noindent In this section we give some examples
symplectic half-flat manifolds and  special Lagrangian
submanifolds. We give also an example of a smooth family of a
symplectic half-flat structures on the $6$-dimensional
torus which is integrable for $t=0$, but not integrable for $t\neq
0$.
\begin{ex}{\em
Let $G$ be the Lie group of matrices of the form
$$
A=\left(
\begin{array}{cccccc}
e^t&0      &xe^t  &0        &0  &y_1\\
0  &e^{-t} &0     &xe^{-t}  &0  &y_2\\
0  &0      &e^t   &0        &0  &w_1\\
0  &0      &0     &e^{-t}   &0  &w_2\\
0  &0      &0     &0        &1  &t\\
0  &0      &0     &0        &0  &1\\
\end{array}
\right)
$$
Let
\begin{equation}
\label{ai}
\begin{aligned}
&\alpha_1=dt,\;\;\alpha_2=dx,\;\;\alpha_3=e^{-t}dy_1-xe^{-t}dw_1\\
&\alpha_4=e^tdy_2-xe^tdw_2,\;\;\alpha_5=e^{-t}dw_1,\;\;\alpha_6=e^tdw_2\,.
\end{aligned}
\end{equation}
Then $\{\alpha_1,\dots,\alpha_6\}$ is a basis of left-invariant
1-forms.\\
By (\ref{ai}) we easily get
\begin{equation}
\begin{cases}
d\alpha_1=d\alpha_2=0\\
d\alpha_3=-\alpha_1\wedge\alpha_3-\alpha_2\wedge\alpha_5\\
d\alpha_4=\alpha_1\wedge \alpha_4-\alpha_2\wedge\alpha_6\\
d\alpha_5=-\alpha_1\wedge\alpha_5\\
d\alpha_6=\alpha_1\wedge\alpha_6\,.
\end{cases}
\end{equation}
Let $\{\xi_1,\dots \xi_6\}$ be the dual frame of
$\{\alpha_1,\dots,\alpha_6\}$; we have
\begin{equation}
\label{xii}
\begin{aligned}
\xi_1=\de{}{t}\,,\;\; \xi_2=\de{}{x}\,,\;\;
\xi_3=e^t\de{}{y_1}\,,\;\;
\xi_4=e^{-t}\de{}{y_2}\\
\xi_5=e^t\de{}{w_1}+xe^t\de{}{y_1}\,,\;\;
\xi_6=e^{-t}\de{}{w_2}+xe^{-t}\de{}{y_2}\,.
\end{aligned}
\end{equation}
From (\ref{xii}) we obtain
\begin{equation}
\label{xii2}
\begin{aligned}
&[\xi_1,\xi_3]=\xi_3,\;\;[\xi_1,\xi_4]=-\xi_4,\;\;,[\xi_1,\xi_5]=\xi_5\\
&[\xi_1,\xi_6]=-\xi_6,\;\;[\xi_2,\xi_5]=\xi_3,\;\;[\xi_2,\xi_6]=\xi_4
\end{aligned}
\end{equation}
and the other brackets are zero.\\
Therefore $G$ is a non-nilpotent solvable Lie group.\\
By \cite{FLS} $G$ has a cocompact lattice  $\Gamma$. Hence
$$
M=G/\Gamma
$$
is a compact solvmanifold of dimension 6. Let us denote with
$\pi\colon \R^6\rightarrow M$ the natural projection. Define
$$
\omega=\alpha_1\wedge\alpha_2+\alpha_3\wedge\alpha_6+\alpha_4\wedge\alpha_5
$$
and
$$
\begin{array}{lll}
J(\xi_1)=\xi_2\,,& J(\xi_3)=\xi_6\,,& J(\xi_4)=\xi_5\\[3pt]
J(\xi_2)=-\xi_1\,,&J(\xi_6)=-\xi_3\,,&J(\xi_5)=-\xi_6\,.
\end{array}
$$
\vskip.2truecm\noindent Then $\omega$ is a symplectic form on $M$
and $J$ is an $\omega$-calibrated almost complex structure on $M$.
Set
$$
\psi=i(\alpha_1+i\alpha_2)\wedge(\alpha_3+i\alpha_6)\wedge(\alpha_4+i\alpha_5)\,;
$$
a direct computation shows that $(\omega,J,\psi)$ is a
symplectic half-flat structure on $M$.\\
Let consider now the lattice $\Sigma\subset\R^4$ given by
$$
\Lambda:=Span_{\Z}\left\{ \left(
\begin{array}{c}
-\mu\\
1\\
0\\
0
\end{array}
\right), \left(
\begin{array}{c}
1\\
\mu\\
0\\
0
\end{array}
\right), \left(
\begin{array}{c}
0\\
0\\
-\mu\\
1
\end{array}
\right), \left(
\begin{array}{c}
0\\
0\\
1\\
\mu
\end{array}
\right) \right\}\,,
$$
where $\mu=\frac{\sqrt{5}-1}{2}$. Let $\T^4$ be the torus
$$
\T^4=\R^4/\Lambda\,.
$$
For any $p,q\in\Z$ let $\rho(p,q)$ be the transformation of $\T^4$
represented by the matrix
$$
\left(
\begin{array}{cccc}
e^{p\lambda}&0            &qe^{p\lambda}&0\\
0           &e^{-p\lambda}&0            &qe^{-p\lambda}\\
0           &0            &e^{p\lambda} &0\\
0           &0            &0            &e^{-p\lambda}

\end{array}
\right)\,,
$$
where $\lambda=\log\frac{3+\sqrt{5}}{2}$.\\
Then
$$
A(p,q)([y_1,y_2,z_1,z_2],(t,x))=(\rho(p,q)[y_1,y_2,z_1,z_2],(t+p,x+q))
$$
is a transformation of $\T^4\times \R^2$ for any $p,q\in\Z$. Let
$\Theta$ be the group of such transformations. The manifold $M$ can
be identified with
\begin{equation}
\label{identification} \frac{\T^4\times\R^2}{\Theta}
\end{equation}
(see \cite{FLS}).\\
Let consider now the involutive distribution ${\mathcal D}$
generated by $\{\xi_2,\xi_3,\xi_4\}$ and
let $p\colon L\hookrightarrow M$ be the leaf through $\pi(0)$.\\
By (\ref{xii}) and the identification (\ref{identification}) we get
$$
\pi^{-1}(L)=\{x=(x_1,\dots,x_6)\in\R^6\;\;|\;\;x_1=x_5=x_6=0\}\,;
$$
hence $L$ is a compact submanifold of $M$. By a direct computation
one can check that
$$
\left\{
\begin{aligned}
& p^*(\omega)=0\,,\\
& p^*(\Im \mathfrak{m}\,\psi)=0\,.
\end{aligned}
\right.
$$
Hence $L$ is a special Lagrangian submanifold of
$(M,\omega,J,\psi)$.  Moreover, by \cite{FLS} we have
$$
H^2(M,\R)=\Span_{\R}\{[\alpha_1\wedge\alpha_2],[\alpha_5\wedge\alpha_6],
[\alpha_3\wedge\alpha_6+\alpha_4\wedge\alpha_5]\}
$$
Therefore, we get
$$
p^*(H^2(M,\R))=0\,.
$$
}
\end{ex}
\begin{ex}
\label{torus}
{\em Let  $(x_1,\dots,x_6)$ be coordinates on $\R^6$ and let
$$
\omega_3=dx_1\wedge dx_4+ dx_2\wedge dx_5+dx_3\wedge dx_6
$$
be the standard symplectic form on $\R^6$. Let $a=a(x_1)$,
$b=b(x_2)$, $c=c(x_3)$ be three smooth functions such that
$$
\lambda_1:=b(x_2)-c(x_3),\quad\lambda_2:=-a(x_1)+c(x_3),\quad\lambda_3=a(x_1)
-b(x_2)
$$
are $\Z^6$-periodic. Let us consider the $\omega_3$-calibrated
complex structure on $\R^6$ defined by
$$
\left\{
\begin{array}{lll}
J(\de{}{x_r})&=&e^{-\lambda_r}\de{}{x_{3+r}}\\[7pt]
J(\de{}{x_{3+r}})&=&-e^{\lambda_r}\de{}{x_r}
\end{array}
\right.
$$
$r=1,2,3$. Define a (3,0)-form on $\R^6$ by
$$
\psi=i(dx_1+ie^{\lambda_1}dx_4)\wedge(dx_2+ie^{\lambda_2}
dx_5)\wedge(dx_3+ie^{\lambda_3}dx_6)\,.
$$
Then we get
$$
\begin{cases}
\psi\wedge\overline{\psi}=-i\frac{4}{3}\omega_3^3\\[5pt]
d\Re \mathfrak{e}\,\psi=0\,.
\end{cases}
$$
Since $\lambda_1,\lambda_2,\lambda_3$ are $\Z^6$-periodic,
$(\omega_3,J,\psi)$ defines a symplectic half-flat
structure on the torus $\T^6=\R^6/\Z^6$. Now consider the
three-torus $L=\pi(X)$, where $\pi\colon\R^6\rightarrow\T^6 $ is the
natural projection and
$$
X=\{(x_1,\dots,x_6)\in \R^6\;\;|\;\;x_1=x_2=x_3=0\}.
$$
It is immediate to check that $L$ is a special Lagrangian
submanifold of $\T^6$.}
\end{ex}
Now we are ready to state the following
\begin{theorem}\label{remtorus}
There exists a family  $(\omega_t,J_t,\psi_t)$ of symplectic half-flat structures on the $6$-dimensional
torus $\T^6$, such that $(\omega_0,J_0,\psi_0)$ is the standard Calabi-Yau structure, but $(\omega_t,J_t,\psi_t)$  is not integrable
for $t\neq 0$.
\end{theorem}
\begin{proof}
By using the notation of example \ref{torus}, let
$$
\left\{
\begin{array}{lll}
J(\de{}{x_r})&=&e^{-t\lambda_r}\de{}{x_{3+r}}\\[7pt]
J(\de{}{x_{3+r}})&=&-e^{t\lambda_r}\de{}{x_r}\,,
\end{array}
\right.
$$
for $r=1,2,3$,
$$
\omega_t=dx_1\wedge dx_4+ dx_2\wedge dx_5+dx_3\wedge dx_6
$$
and
$$
\psi_t=i(dx_1+ie^{t\lambda_1}dx_4)\wedge(dx_3+ie^{t\lambda_2}
dx_5)\wedge(dx_3+ie^{t\lambda_3}dx_6)\,.
$$
Then $(\T^6,\omega_t,J_t,\psi_t)$ is a symplectic
half-flat manifold for any $t\in\R$, such that
$(\T^6,\omega_0,J_0,\psi_0)$ is the standard holomorphic
Calabi-Yau torus and $J_t$ is non-integrable for $t\neq 0$ (here we
assume that $\lambda_1,\lambda_2,\lambda_3$ are not constant).
\end{proof}
\begin{ex}{\em
Let consider now the Lie group $G$ of matrices of the form
$$
A=\left(
\begin{array}{cccccc}
1 &0 &x_1  &u_1  &0  &0\\
0 &1 &x_2  &u_2  &0  &0\\
0 &0 &1    &y    &0  &0\\
0 &0 &0    &1    &0  &0\\
0 &0 &0    &0    &1  &t\\
0 &0 &0    &0    &0  &1\\
\end{array}
\right)
$$
where $x_1,x_2,u_1,u_2,y,t$ are real numbers. Let $\Gamma$ be the
subgroup $G$ formed by the matrices having integral entries. Since
$\Gamma$ is a cocompact lattice of $G$, then
$M:=G/\Gamma$ is a $6$-dimensional nilmanifold.\\
Let consider
\begin{equation*}
\begin{aligned}
&\xi_1=\de{}{y}+x_1\de{}{u_1}+x_2\de{}{u_2}\,,\;\;\xi_2=
\de{}{x_2}\,,\\
&\xi_3=\de{}{x_1}\,,\;\;
\xi_4=\de{}{t}\,,\;\;\xi_5=\de{}{u_1}\,,\;\;\xi_6=\de{}{u_2}\,.
\end{aligned}
\end{equation*}
Then $\{\xi_1,\dots,\xi_6\}$ is a $G$-invariant global frame on
$M$.\\
The respective coframe $\{\alpha_1,\dots,\alpha_6\}$ satisfies
\begin{equation}\label{alpha3}
\begin{cases}
d\alpha_1=d\alpha_2=d\alpha_3=d\alpha_4=0\\
d\alpha_5=\alpha_1\wedge\alpha_3\\
d\alpha_6=\alpha_1\wedge\alpha_2\,.
\end{cases}
\end{equation}
The symplectic half-flat structure on $M$ is given by the
symplectic form
$$
\omega=\alpha_1\wedge\alpha_4+\alpha_2\wedge\alpha_5+\alpha_3\wedge\alpha_6\,,
$$
by the $\omega$-calibrated almost complex structure
$$
\begin{array}{lll}
J(\xi_1)=\xi_4\,,&J(\xi_2)=\xi_5\,,&J(\xi_3)=\xi_6\,,\\[3pt]
J(\xi_4)=-\xi_1\,,&J(\xi_5)=-\xi_2\,,&J(\xi_6)=-\xi_3
\end{array}
$$
and by the complex volume form
$$
\psi=(\alpha_1+i\alpha_4)\wedge(\alpha_2+i\alpha_5)\wedge(\alpha_3+i\alpha_6)\,.
$$
By a direct computation we get
$$
\begin{aligned}
&\Re \mathfrak{e}\,\psi=\alpha_{123}-\alpha_{345}+\alpha_{246}-
\alpha_{156}\,,\\
&\Im \mathfrak{m}\,\psi=\alpha_{234}-\alpha_{135}+\alpha_{126}
-\alpha_{456};\,.
\end{aligned}
$$
Let
$$
X=\{A\in G\;\;\vert\;\;y=x_2=u_2=0\}
$$
and
$$
L=\pi(X)\,,
$$
$\pi\colon G\rightarrow M$ being the canonical projection. Then $L$
is a special Lagrangian torus embedded in $(M,\omega,J,\psi)$.}
\end{ex}
In order to obtain some cohomological obstructions to the existence of a symplectic half-flat structure $(\omega,J,\psi)$
on a compact $6$-manifold $M$, one can ask if the cohomology class $[\Re\mathfrak{e}\,\psi]$ is always non-trivial. Observe that in the
Calabi-Yau case one has $[\Re\mathfrak{e}\,\psi]\neq 0$. In our context we have the following
\begin{theorem}\label{principale}
There exists a compact $6$-dimensional manifold admitting a symplectic half-flat structure $(\omega,J,\psi)$ such that
$$
[\Re\mathfrak{e}\,\psi]=0\,.
$$
\end{theorem}
\begin{proof}
Let $G$ be the Lie group consisting of matrices of the form
$$
A=\left(
\begin{array}{cccc}
e^{\lambda z}&0            &0            &x\\
0           &e^{-\lambda z}&0            &y\\
0           &0             &1            &z\\
0           &0             &0            &1
\end{array}
\right)\,,
$$
where $x,y,z$ are real numbers and
$$
\lambda=\log\frac{3+\sqrt{5}}{2}\,.
$$
Then $G$ is a connected solvable Lie group admitting a cocompact lattice $\Gamma$ (see \cite{FG}). Therefore $N=G/\Gamma$
is a 3-dimensional parallelizable solvmanifold. It can be easily showed  (see \cite{FG} again) that
there exists a coframe $\{\alpha_1,\alpha_2,\alpha_3\}$ on $M$ satisfying the following
structure equations
$$
\begin{cases}
d\alpha_1=-\lambda\alpha_{1}\wedge\alpha_3\\
d\alpha_2=\lambda\alpha_2\wedge\alpha_3\\
d\alpha_3=0\,.
\end{cases}
$$
Let $M=N\times N$. Then $M$ is a compact $6$-manifold admitting a coframe $\{\alpha_1,\dots,\alpha_{6}\}$ satisfying
$$
\begin{cases}
d\alpha_1=-\lambda\alpha_{1}\wedge\alpha_3\\
d\alpha_2=\lambda\alpha_{2}\wedge\alpha_3\\
d\alpha_3=0\\
d\alpha_4=-\lambda\alpha_{4}\wedge\alpha_6\\
d\alpha_5=\lambda\alpha_5\wedge\alpha_6\\
d\alpha_6=0\,.
\end{cases}
$$
Let us now construct a symplectic half-flat structure on $M$ satisfying $[\Re\mathfrak{e}\,\psi]=0$.
Let $(\omega,J)$ be the almost K\"ahler structure on $M$ given by the symplectic form
$$
\omega=\alpha_1\wedge\alpha_2+\alpha_4\wedge\alpha_5+\alpha_3\wedge\alpha_6
$$
and by the $\omega$-calibrated
almost complex structure $J$ defined as
$$
\begin{array}{lll}
J(\xi_1)=\xi_2\,,& J(\xi_3)=\xi_6\,,& J(\xi_4)=\xi_5\\[3pt]
J(\xi_2)=-\xi_1\,,&J(\xi_6)=-\xi_3\,,&J(\xi_5)=-\xi_6\,,
\end{array}
$$
where $\{\xi_1,\dots,\xi_6\}$ is the frame on $M$ dual to
$\{\alpha_1,\dots,\alpha_6\}$.\\
Then the complex 3-form
$$
\psi=\frac{\sqrt{2}}{2}(1-i)\,(\alpha_1+i\alpha_2)\wedge(\alpha_4+i\alpha_5)\wedge(\alpha_3+i\alpha_6)\,.
$$
defines together with $(\omega,J)$ a symplectic half-flat structure on $M$.\\
Moreover
$$
\begin{aligned}
\Re\mathfrak{e}\,\psi=\frac{\sqrt{2}}{2}(-\alpha_{134}+\alpha_{146}-\alpha_{135}-
\alpha_{156}-\alpha_{234}-\alpha_{246}+
\alpha_{235}-
\alpha_{256})
\end{aligned}
$$
and a direct computation gives
$$
\Re\mathfrak{e}\,\psi=\frac{\sqrt{2}}{2\lambda}\,d(\alpha_1\wedge\alpha_4+\alpha_1\wedge\alpha_5-\alpha_2\wedge\alpha_4+\alpha_2\wedge\alpha_5)\,.
$$
This ends the proof.
\end{proof}
\begin{rem}{\rm
The symplectic manifold $(M,\omega)$, described in the last example, satisfies the {\em hard Lefschetz condition}, i.e.
$$
\begin{array}{lclc}
\omega^k :&\Lambda^{3-k}(M)&\to      &\Lambda^{3+k}(M)\\
          &\alpha          &\mapsto  &\omega^k\wedge\alpha
\end{array}
$$
$k=1,2$ induces an isomorphism in cohomology. Indeed, it is immediate
to check that
\begin{eqnarray*}
H^1(M,\R)&\!\!\!=\!\!\!&\Span_\R\{[\alpha_3],[\alpha_6]\}\,,\\[3pt]
H^2(M,\R)&\!\!\!=\!\!\!&\Span_\R\{[\alpha_1\wedge\alpha_2],[\alpha_4\wedge\alpha_5],[\alpha_3\wedge\alpha_6]\}\,,\\[3pt]
H^4(M,\R)&\!\!\!=\!\!\!&\Span_\R\{[\alpha_3\wedge\alpha_4\wedge\alpha_5\wedge\alpha_6],
[\alpha_1\wedge\alpha_2\wedge\alpha_3\wedge\alpha_6],
[\alpha_1\wedge\alpha_2\wedge\alpha_4\wedge\alpha_5]\}\,,\\[3pt]
H^5(M,\R)&\!\!\!=\!\!\!&\Span_\R\{[\alpha_1\wedge\alpha_2\wedge\alpha_4\wedge\alpha_5\wedge\alpha_6],
[\alpha_1\wedge\alpha_2\wedge\alpha_4\wedge\alpha_5\wedge\alpha_6]\}\,.
\end{eqnarray*}
Therefore, if $\omega =\alpha_1\wedge\alpha_2+ \alpha_4\wedge\alpha_5+\alpha_3\wedge\alpha_6$, then it is easy to
verify that $(M,\omega)$ satisfies the hard Lefschetz condition. Hence, in view of a Proposition of Hitchin
(see \cite[prop. 7]{Hi3}) we have that $(M,\omega)$ satisfies a
$dd^J$-lemma (see \cite[def. 5]{Hi3} for the precise definition). Consequently, the manifold $M$ has a symplectic half-flat structure $(\omega,J,\psi)$ such that
the symplectic structure $\omega$ gives rise to a generalized Calabi-Yau structure on $M$ satisfying the $dd^J$-lemma.}
\end{rem}
\section{The four-dimensional case}
Let $(M,\omega)$ be a (compact) $4$-dimensional symplectic manifold and $J$ be an
$\omega$-calibrated almost complex structure on $M$. Let $\psi$
be a nowhere vanishing $(2,0)$-form on $M$ satisfying
$$
\psi\wedge\overline{\psi}=2\,\omega^2\,.
$$

Then the conditions
$$
\begin{cases}
\nabla\psi=0\\
d\Re \mathfrak{e}\,\psi=0
\end{cases}
$$
imply
$$
d(\Im \mathfrak{m} \,\psi)=0\,.
$$
Indeed, if $\nabla\psi =0$, then $\overline{\partial}_J\psi =0$ and
$$
d\Re \mathfrak{e}\,\psi=0 \Longrightarrow  \overline{\partial}_J\psi +\partial_J\overline{\psi} +
\overline{A}_J\psi + A_J\overline{\psi}=\overline{A}_J\psi + A_J\overline{\psi}=0\,.
$$
Since $\overline{A}_J\psi \in\Lambda_J^{1,2}(M)$ and $A_J\overline{\psi} \in\Lambda_J^{2,1}(M)$ , we get
$d\psi =0$ which implies that $J$ is holomorphic.
In dimension 4 we adopt the following definition (see also \cite{dB})
\begin{definition}
Let $M$ be a {\rm (}compact{\rm )} $4$-dimensional manifold. A {\em special symplectic structure} on $M$
is a triple $(\omega,J,\psi)$, where
\begin{itemize}
\item $\omega$ is a symplectic form,\vskip.2truecm\noindent
\item $J$ is an $\omega$-calibrated almost complex structure on $M$,\vskip.2truecm\noindent
\item $\psi$ is a non-vanishing $(2,0)$-form satisfying\vskip.2truecm\noindent
$$
\begin{cases}
\psi\wedge\overline{\psi}=2\,\omega^2\\
d\Re \mathfrak{e}\,\psi=0\,.
\end{cases}
$$
\end{itemize}
\end{definition}
The following proposition gives a characterizations of a
$4$-dimensional special symplectic structure in terms of
differential forms. The proof of this proposition can be obtained
using the same argument as in \cite[proposition 1]{ContiS}.
\begin{prop}
\label{SU2}
Symplectic half-flat structures on a $4$-dimensional manifold are
in one-to-one correspondence to the triple
$(\omega,\Omega_+,\Omega_-)$ of $2$-forms on $M$
satisfying the following properties
\begin{enumerate}
\item[1.]$\omega\wedge\Omega_+=\omega\wedge\Omega_-=\Omega_+\wedge\Omega_-=0\,;$
\vskip0.2cm
\item[2.]$\Omega_+\wedge\Omega_+=\Omega_-\wedge\Omega_-=\omega\wedge\omega\neq 0\,;$
\vskip0.2cm
\item[3.]if $\iota_X\Omega_+=\iota_Y\Omega_-$, then $\omega(X,Y)\geq 0\,;$
\vskip0.2cm
\item[4.] $d\omega=d\Omega_+=0\,.$
\end{enumerate}
\end{prop}
Note that, when a triple $(\omega,\Omega_+,\Omega_-)$ is given, the complex volume form of the associated special
symplectic structure is
simply obtained by taking
$$
\psi=\Omega_++i\Omega_-\,.
$$
The following proposition gives an explicit formula for the almost complex structure induced by a triple of $2$-forms
$(\omega,\Omega_+,\Omega_-)$ satisfying the properties stated above
\begin{prop}
Let $M$ be a $4$-manifold and let $(\omega,\Omega_+,\Omega_-)$ be a triple of forms satisfying the properties $1-4$ of proposition $\ref{SU2}$. Then
the $\omega$-calibrated almost complex structure induced
by $(\omega,\Omega_+,\Omega_-)$ is the dual of the endomorphism $P$ of $T^*M$ defined by the
following formula
$$
P(\alpha)=\bigstar(\omega\wedge\bigstar(\Omega_+\wedge\alpha))\,,
$$
where, as usual, $\bigstar$ is the symplectic star operator associated with $\omega$ and $\alpha\in T^*M$.
\end{prop}
\begin{proof}
Since the triple $(\omega,\Omega_+,\Omega_-)$ induces an SU(2)-structure, $M$ admits
a local coframe $\mathcal{B}=\{e^1,e^2,e^3,e^4\}$ with respect to which the structure forms $\omega,\Omega_+,\Omega_-$ take the standard expressions
$$
\begin{aligned}
\omega & =e^1\wedge e^2+e^3\wedge e^4\,,\\[3pt]
\Omega_+ & =e^1\wedge e^3-e^2\wedge e^4\,,\\[3pt]
\Omega_-& =e^1\wedge e^4+e^2\wedge e^3\,.
\end{aligned}
$$
With respect to the dual frame of $\mathcal{B}$ the almost complex structure $J$ induced by $(\omega,\Omega_+,\Omega_-)$ on $M$
reads as the standard one
given by the following relations
$$
\begin{aligned}
&J(e_1)=e_2\,,\quad J(e_2)=-e_1\\
&J(e_3)=e_4\,,\quad J(e_4)=-e_3\,.
\end{aligned}
$$
Hence it is sufficient to prove the proposition for these data and the latter is a direct computation.
\end{proof}
Note that $\Omega_+,\Omega_-$ and $J$ are related by $\Omega_-=J\Omega_+$. Furthermore if $*$ denotes the Hodge star operator of the metric
induced by a triple $(\omega,\Omega_+,\Omega_-)$, then $*\Omega_+=\Omega_-$.
The following lemma gives a topological obstruction to the existence of special symplectic structures on compact $4$-manifolds.
\begin{lemma}
Let $M$ be a $4$-dimensional compact
manifold admitting a special symplectic structure. Then
$$
\dim(H^2(M,\R))\geq 2\,.
$$
\end{lemma}
\begin{proof}
Let $(\omega,J,\psi)$ be a special symplectic structure on $M$ and let $\Omega_+=\Re\mathfrak{e}\,\psi$,
$\Omega_-=\Im\mathfrak{m}\,\psi$. First of all we observe that $\Omega_+$ is a symplectic form on $M$ and consequently it cannot be exact.
Furthermore if $a[\omega]+b[\Omega_+]=0$ for some $a,b\in \R$, then
$$
a\Omega_++b\omega=d\alpha\,,
$$
for some $\alpha\in\Lambda^1(M)$. This last equation together with $\Omega_+\wedge\omega=0$ readily implies $b\omega^2=d(\alpha\wedge\omega)$, which
forces $b$ to vanish. Hence $\omega$ and $\Omega_+$ induce $\R$-linear independent classes in $H^2(M,\R)$ and $\dim(H^2(M,\R))\geq 2$.
\end{proof}
As in the $6$-dimensional case, if $(M,\omega,J,\psi)$ is a
$4$-dimensional symplectic half-flat manifold, then
$\Omega=\Re\mathfrak{e}\,\psi$ is a calibration on $M$.
Furthermore in an obvious way we have the definition of special
Lagrangian submanifold.\\
We end this section with the following
\begin{ex}{\em
Now we recall the construction of the Kodaira-Thurston manifold.\\
Let $G$ be the Lie subgroup of ${\rm GL}(5,\R)$ whose matrices have
the following form
$$
A=\left(
\begin{array}{cccccc}
1&x&z&0&0\\
0&1&y&0&0\\
0&0&1&0&0\\
0&0&0&1&t\\
0&0&0&0&1
\end{array}
\right)
$$
where $x,y,z,t\in\R$. Let $\Gamma$ be the subgroup of $G$ of
matrices with integers entries. Since $\Gamma$ is a cocompact
lattice in $G$, we get that
$$
M=G/\Gamma
$$
is a compact manifold.
$M$ is called the {\em Kodaira-Thurston manifold}.\\
Let $\{\xi_1,\dots,\xi_4\}$ be the global frame of $M$ given by
$$
\xi_1=\de{}{x}\,,\;\;\xi_2=\de{}{y}+x\de{}{z}\,,\;\;\xi_3=\de{}{z}\,,\;\;\xi_4=\de{}{t}\,.
$$
We easily get
$$
[\xi_1,\xi_2]=\xi_3
$$
and the other brackets are zero.\\
The dual frame of $\{\xi_1,\dots,\xi_4\}$ is given by
$$
\alpha_1=dx\,,\;\;\alpha_2=dy\,,\;\;\alpha_3=dz-x
dy\,,\;\;\alpha_4=dt\,.
$$
We have
$$
d\alpha_1=d\alpha_2=d\alpha_4=0\,,\;\;d\alpha_3=-\alpha_1\wedge\alpha_2\,.
$$
The special symplectic structure on $M$ is given by the forms
$$
\begin{cases}
\omega=\alpha_1\wedge\alpha_3+\alpha_2\wedge\alpha_4\\
\psi=i(\alpha_1+i\alpha_3)\wedge(\alpha_2+i\alpha_4)
\end{cases}
$$
and by the almost complex structure
$$
\begin{aligned}
&J(\xi_1)=\xi_3\,,\;\;&J(\xi_2)=\xi_4\,,\\
&J(\xi_3)=-\xi_1\,,\;\;      &J(\xi_4)=-\xi_2\,.
\end{aligned}
$$
We immediately get
$$
\begin{aligned}
\Im \mathfrak{m}\, \psi=\alpha_1\wedge\alpha_2-\alpha_3\wedge\alpha_4\,,\\
\Re \mathfrak{e}\,
\psi=\alpha_2\wedge\alpha_3-\alpha_1\wedge\alpha_4\,.
\end{aligned}
$$
Hence
$$
\begin{aligned}
&d\Re \mathfrak{e}\,\psi=0\,.
\end{aligned}
$$
Let $X\subset G$ be the set
$$
X=\{A\in G\;\;\vert\;\;x=t=0\}
$$
and
$$
L=\pi(X),
$$
where $\pi:G\rightarrow M$ is the natural projection. Hence $L$ is a
compact manifold embedded in $M$. Moreover the tangent bundle to $L$
is generated by $\{\xi_2,\xi_3\}$; so we get
$$
\begin{aligned}
&p^*(\omega)=0\,,\\
&p^*(\Im \mathfrak{m}\,\psi)=0\,.
\end{aligned}
$$
Hence $L$ is a special Lagrangian torus.}
\end{ex}
\vskip.5cm \noindent \textbf{Acknowledgements.} This paper has
originated by a discussion with Gang Tian during the Differential
Geometry semester held at Centro De Giorgi Pisa 2004. We are
grateful to him. \\
We also would like to thank Paolo de Bartolomeis for useful remarks. Finally we thank Tommaso Pacini for useful suggestions to
improve the present work.


\begin{thebibliography}{12}
\bibitem{CF}
Chiossi S., Fino A.: Conformally parallel $G_2$ structures on a class of solvmanifolds.  \emph{Math. Z.}  {\bf 252} (2006),
pp. 825--848.
\bibitem{CS} Chiossi S., Salamon S.: The intrinsic torsion of
$\rm SU(3)$ and $G\sb 2$ structures, {\em Differential geometry}, Valencia, 2001,
World Sci. Publishing, River Edge, NJ (2002) pp. 115--133.
\bibitem{CT}
Conti D., Tomassini A.: Special Symplectic $6$-manifolds, to appear in \emph{Q. J. Math.},
{\tt
e-print math.DG/0601002} (2006).
\bibitem{ContiS}
Conti D., Salamon S.: Generalized Killing spinors in dimension 5,
e-print math.DG/0508375, to appear in \emph{Trans. Amer. Math. Soc.}
\bibitem{dB} de Bartolomeis P.: Geometric Structures on Moduli Spaces of Special Lagrangian Submanifolds,
\emph{Ann. di Mat. Pura ed Applicata}, IV, Vol. CLXXIX, (2001),
pp. 361--382.
\bibitem{dBT2} de Bartolomeis P., Tomassini A.:
On the Maslov Index of Lagrangian Submanifolds of Generalized
Calabi-Yau Manifolds, {\em Int. J. of Math.} {\bf 17}, (2006) pp. 921--947.
\bibitem{dBT3} de Bartolomeis P., Tomassini A.:
On solvable Generalized Calabi-Yau Manifolds, {\em Ann. Inst.
Fourier} {\bf 56}, (2006) pp. 1281--1296.
\bibitem{FLS} Fern\'andez M., de Le\'on M., Saralegui M.: A six dimensional Compact Symplectic Solvmanifold without
K\"ahler Structures, {\em Osaka J. Math} {\bf 33} (1996) pp.
19--34.
\bibitem{FG} Fern\'andez M., Gray A.: Compact symplectic sovmanifold not admitting complex structures, {\em Geometria Dedicata} {\bf 34}
(1990), pp. 295--299.
\bibitem{Gau}
Gauduchon P.: Hermitian connections and Dirac operators.
\emph{Boll. Un. Mat. Ital. B} (7) {\bf 11} (1997), no. 2, suppl., pp. 257--288.
\bibitem{G} Gualtieri M.: Generalized complex geometry, DPhil thesis, University of Oxford, 2003, {\tt e-print math.DG/0401221}.
\bibitem{HL} Harvey R., Lawson H. Blaine Jr.:
Calibrated geometries, {\em Acta Math.} {\bf 148} (1982) pp. 47--157.
\bibitem{Hi1} Hitchin N.J.: {\em Stable forms and special metrics},
Global differential geometry: the mathematical legacy of Alfred
Gray (Bilbao, 2000), 70--89, Contemp. Math., 288, Amer. Math.
Soc., Providence, RI, 2001.
\bibitem{Hi3} Hitchin N.J.: Generalized Calabi-Yau Manifolds, {\em Quart. J. Math.} {\bf 54} (2003) pp.
281--308.
\end{thebibliography}
\end{document}